\documentclass[10pt]{article}
\usepackage{calc}
\usepackage[all]{xy}
\usepackage[centertags]{amsmath}
\usepackage{latexsym}
\usepackage{amsfonts}
\usepackage{amssymb}
\usepackage{amsthm}
\usepackage{fancyhdr}
\usepackage{url}
\usepackage{times}
\usepackage{color}
\usepackage{amsmath}
\usepackage{graphicx}
\usepackage[pdfpagelabels,plainpages=false]{hyperref}

\newtheorem{thm}{Theorem}[section]

\newtheorem{prop}[thm]{Proposition}
\newtheorem{lem}[thm]{Lemma}

\allowdisplaybreaks
\hypersetup{ linkcolor =blue,
 citecolor = cyan,
 urlcolor = blue,
colorlinks = true, bookmarksopen = true, bookmarksnumbered = true,
bookmarksopenlevel = {1}, linktocpage = true,
pdfstartview = FitH,
pdftitle = {Linear Case},
pdfsubject
= {article},
pdfauthor = {Habiba  KNANI},
pdfcreator = {pdftextify}, pdfproducer = {pdftextify}
}
\begin{document}

\title{Linear Backward Stochastic Differential Equations with Gaussian
Volterra processes}
\author{H.Knani \thanks{%
Laboratoire de Math\'{e}matiques: Mod\'elisation D\'{e}terministe et Al\'{e}atoire, Hammam
Sousse, Tunisie et UMR-CNRS 7502, Institut \'{E}lie Cartan de Lorraine,
Nancy, France}, M. Dozzi \thanks{%
UMR-CNRS 7502, Institut \'{E}lie Cartan de Lorraine, Nancy, France}}
\maketitle

\begin{abstract}
{\footnotesize Explicit solutions for a class of linear backward stochastic differential equations (BSDE) driven by Gaussian Volterra processes are given. These processes include the multifractional brownian motion and the multifractional Ornstein-Uhlenbeck process. By an It\^o formula, proven in the context of Malliavin calculus, the BSDE is associated to a linear second order partial differential equation with terminal condition whose solution is given by a Feynman-Kac type formula. An application to self-financing trading strategies is discussed.} \newline


\end{abstract}

\section{Introduction}

\quad A\ backward stochastic differential equation (BSDE) with a generator 
$f:[0,T]\times \mathbb{R}\times\mathbb{R}^n\rightarrow\mathbb{R},$ a
terminal value $\xi _{T}$ and driven by a stochastic process $%
X=(X^{1},...,X^{n})$ is given by the equation%
\begin{equation}
Y_{t}=\xi _{T}-\int_{t}^{T}f(s,Y_{s},Z_{s})ds+\int_{t}^{T}Z_{s}dX_{s},\
0\leqslant t\leqslant T.  \label{EDSR}
\end{equation}

\quad A solution is a pair of square integrable processes $Y$ and $%
Z=(Z^{1},...,Z^{n})$ that are adapted to the filtration generated by $X.$
Such equations appear especially in the context of asset pricing and hedging
theory in finance and in the context of stochastic control problems. BSDEs
may be considered as an alternative to the more familiar partial
differential equations (PDE) since the solutions of BSDEs are closely
related to classical or viscosity solutions of associated PDEs. As a
consequence, BSDEs may be used for the numerical solution of nonlinear PDEs.

\quad BSDEs driven by brownian motion have been studied extensively after the
first general existence and uniqueness result proved by E. Pardoux and S.G.\
Peng \cite{PaPe}. For a synthesis of this research work we may refer to the
recent textbooks \cite{Carmona}, \cite{Crepey}, \cite{Pardoux.Rasc}, \cite%
{Pham}, \cite{Touzi}, \cite{ZhangI}. More recently BSDEs driven by
fractional brownian motions have been investigated (see, e.g., \cite%
{Bender05}, \cite{HuPeng},\ \cite{HuOcSo}, \cite{JaBo}, \cite{MatNie}, \cite%
{SoDi}, \cite{WeSh}, \cite{ZhangII}). Since fractional brownian motions
(with Hurst index $H\in (0,1/2)\cup (1/2,1))$ are neither martingales nor
Markov processes, new methods have been developped to show the wellposedness
of BSDEs in certain function spaces. Hereby the stochastic integral $%
\int_{t}^{T}Z_{s}dX_{s}$ has frequently been defined as a Skorohod integral
in the framework of Malliavin calculus, and the notion of quasi-conditional
expectation has been introduced, since the classical notion of conditional
expectation seems not to be convenient for a proof of the existence and
uniqueness of solutions to BSDEs whose driving process is not a martingale (%
\cite{Nualart}, Lemma 1.2.5). Very few articles are concerned with BSDEs for
more general gaussian processes (\cite{B2014}, \cite{BeVi}) or in the
context of the theory of rough paths \cite{DiFri}. In \cite{B2014} the
stochastic integral $\int_{t}^{T}Z_{s}dX_{s}$ is understood in the Wick-It%
\^{o} sense, and the existence and uniqueness of the solution of (\ref{EDSR}%
) is proved for a class of gaussian processes which includes fractional
Brownian motion. The proof is based on a transfer theorem that aims to
reduce the question of wellposedness to BSDEs driven by brownian motion. In \cite%
{BeVi} it is shown that the wellposedness of linear BSDEs with general
square integrable terminal condition $\xi $ holds true if and only if $X$ is
a martingale.

\quad This paper is concerned with linear BSDEs driven by gaussian Volterra
processes $X$. This class of processes contains multifractional brownian
motions and multifractional Ornstein-Uhlenbeck processes. Contrary to
fractional Brownian motion where the Hurst parameter $H$ is constant, it
becomes for multifractional Brownian motion a function $h$ which is assumed
here to be differentiable and with values in $(1/2,1)$. The aim is to obtain
the solution of the linear BSDE with the associated linear PDE whose
classical solution is given explicitely. This generalizes a result in 
\cite{Bender05} obtained for fractional brownian motion. We define the
stochastic integral $\int_{t}^{T}Z_{s}dX_{s}$ as a divergence integral, and
extend an It\^{o} formula in \cite{AMN} to the multidimensional case. The It%
\^{o} formula is then applied to the solution of the associated PDE in order
to get a solution of the BSDE. The question of uniqueness of solutions of
(not necessarily linear) BSDEs driven by Volterra processes will be treated
in a separate paper. Special attention is given to the fact that the
variance of Volterra processes is not necessarily an increasing function of
time, but in general only of bounded variation. The explicit solution of the
associated PDE contains this variance and is given by a Feynman-Kac type
formula on time intervals where it is increasing. The application of this
formula to the BSDEs is therefore restricted to time intervals where this variance is increasing.

\quad In this section we define the class of Volterra processes $X$ we have
in mind and the linear BSDEs and the associated PDE. Section 2 is concerned
with complements on the Skorohod integral with respect to Volterra
processes. The It\^{o} formula is proved in Section 3 and applied in Section
4 to the linear BSDE.

\subsection{Gaussian Volterra processes}

Let $X=\{X_{t},0\leqslant t\leqslant T\}$ be a zero mean continuous Gaussian
process given by 
\begin{equation}
X_{t}=\displaystyle\int_{0}^{T}K(t,s)dW_{s}  \label{X}
\end{equation}%
where $W=\{W_{t},0\leqslant t\leqslant T\}$ is a standard Brownian motion
and $K:[0,T]^{2}\rightarrow 
\mathbb{R}
$ is a square integrable kernel, i.e. $\int_{[0,T]^{2}}K(t,s)^{2}dtds<+%
\infty .$ We assume that $K$ is of Volterra type, i.e, $K(t,s)=0$ whenever $%
t<s$. Usually, the representation (\ref{X}) is called a Volterra
representation of $X$. Gaussian Volterra processes and their stochastic
analysis have been studied e.g. in \cite{AMN}, \cite{Coup.Maslk} and \cite%
{SoVi}. In \cite{AMN} $K$ is called \textit{regular }if it satisfies\textit{%
\ }

\textit{(H) } For all $s\in (0,T],$ $\int_{0}^{T}\mid K\mid
((s,T],s)^{2}ds<\infty ,$ where $\mid K\mid ((s,T],s)$ denotes the total
variation of $K(.,s)$ on $(s,T]$.

We assume the following condition on $K(t,s)$ which is more restrictive than 
\textit{(H) }(\cite{AMN}, \cite{Coup.Maslk}):

\emph{(H1) \ }$K(t, s)$ is continuous for all $ 0<s \leqslant t<T $ and  continuously differentiable in the variable $t$ in 
${0<s< t<T}$,

\emph{(H2) \ }For some $0<\alpha, \beta <\frac{1}{2},$ there is a finite
constant $c>0$ such that
\begin{equation}
\begin{aligned} \Big| \frac{\partial K}{\partial t}(t, s)\Big| \leqslant
c(t-s)^{\alpha-1}\Big(\frac{t}{s}\Big)^{\beta},\,\text{for all $0<s<t<T.$}
\end{aligned}  \notag
\end{equation}

\quad The covariance function of $X$ is given by 
\begin{equation}
R(t,s):=\mathbb{E}X_{t}X_{s}=\displaystyle\int_{0}^{inf(t,s)}K(t,u)K(s,u)du.
\label{R}
\end{equation}%
We discuss shortly some examples of Gaussian Volterra processes that
satisfy\ \emph{(H1) and (H2)}.

\textit{Example 1.} The multi-fractional brownian motion (mbm) ($%
B_{t}^{h(t)},0\leqslant t\leqslant T)$ with Hurst function $%
h:[0,T]\rightarrow \lbrack a,b]\subset (\frac{1}{2},1)$. Its kernel is given
by \cite{Bouf.Dozzi} 
\begin{equation}
K(t,s)=s^{1/2-h(t)}\displaystyle\int_{s}^{t}(y-s)^{h(t)-3/2}y^{h(t)-1/2}dy,
\label{Kmbm}
\end{equation}

where $h$ is assumed to be continously differentiable with bounded
derivative. We get 
\begin{equation*}
\frac{\partial K}{\partial t}(t,s)=h^{\prime }(t)s^{\frac{1}{2}%
-h(t)}\int_{s}^{t}(y-s)^{h(t)-\frac{3}{2}}y^{h(t)-\frac{1}{2}}ln((\frac{y}{s}%
-1)y)dy+s^{\frac{1}{2}-h(t)}(t-s)^{h(t)-\frac{3}{2}}t^{h(t)-\frac{1}{2}}.
\end{equation*}
A straightforward calculation shows that (H2) is satisfied with $\alpha=a- 
\frac{1}{2},$ $\beta=b+\epsilon-\frac{1}{2}$ with $\epsilon$ small enough
and $c$ depends on $a$, $b$, $T$ and $\epsilon.$ The mbm generalizes
fractional brownian motion (fbm) with Hurst index $H>1/2.$ Mbm is a more
flexible model than fbm since the H\"{o}lder continuity of its trajectories
varies with $h.$ The trajectories of mbm $B_{\cdot }^{h(\cdot )}$ are in
fact locally H\"{o}lder continous of order $h(t)$ at $t.$

\textit{Example 2.} The multi-fractional Ornstein-Uhlenbeck process $%
U=\{U_{t},0\leqslant t\leqslant T\}$ given by $U_{t}=\displaystyle%
\int_{0}^{t}e^{-\theta (t-s)}dB_{s}^{h(s)},$ where $\theta >0$ is a
parameter and $B^{h}$ is the mbm of Example 1$.$ The kernel of $U$ is given
by \ 
\begin{equation}
\mathcal{K}(t,r)=\displaystyle\int_{r}^{t}e^{-\theta (t-s)}\frac{\partial K}{%
\partial s}(s,r)ds=K(t,r)-\theta \displaystyle\int_{r}^{t}e^{-\theta
(t-s)}K(s,r)ds.  \notag
\end{equation}

In fact we have%
\begin{eqnarray*}
\int_{0}^{t}\mathcal{K}(t,r)dW_{r} &=&\displaystyle\int_{0}^{t}K(t,r)dW_{r}-%
\theta \displaystyle\int_{0}^{t}\left( \int_{r}^{t}e^{-\theta
(t-s)}K(s,r)ds\right) dW_{r} \\
&=&B_{t}^{h(t)}-\theta \displaystyle\int_{0}^{t}e^{-\theta
(t-s)}B_{s}^{h(s)}ds.
\end{eqnarray*}

An integration by parts gives the representation of $U.$ \ We notice that in
the framework of the divergence integral (Section 2) the integral with
respect to mbm can be reduced to an integral with respect to brownian
motion. (H2) is satisfied with the same values of $\alpha$ and $\beta$ as in
Example 1. \emph{\ }

\textit{Example 3. }The Liouville multi-fractional Brownian motion (LmBm) $%
(B_{t}^{L,h(t)},t\in \lbrack 0,T])$ with Hurst function $h$ as in Example 1.
Its kernel is given by $\tilde{\mathcal{K}}(t,r)=(t-r)^{h(t)-\frac{1}{2}%
}1_{(0,t]}(r).$

\subsection{Linear backward stochastic differential equations}

Let $W=(W^{1},...,W^{n})$ a standard brownian motion in $%
\mathbb{R}
^{n},$ defined on the probability space $(\Omega ,\mathcal{F},P),$ and let $%
\mathbb{F=\{}\mathcal{F}_{t}\subset \mathcal{F},$ $t\in \lbrack 0,T]\}$ be
the filtration generated by $W$ and augmented by the $P$-null sets$.$ We
consider the $%
\mathbb{R}
^{n}$-valued Volterra processes $X=(X^{1},...,X^{n})$ given by 
\begin{equation}
X_{t}^{j}=\displaystyle\int_{0}^{t}K^{j}(t,s)dW_{s}^{j},\ \text{$j=1,...,n,$}
\label{Xj}
\end{equation}%
where $K^{j}:[0,T]^{2}\rightarrow 
\mathbb{R}
$ satisfies the conditions\emph{\ (H1)} and\emph{\ (H2)}. Let $\sigma ^{j},$ 
$j=1,...,n$ be bounded functions on $[0,T],$ and let $b^{j}\in \mathcal{C}%
^{1}((0,T),%
\mathbb{R}
)\cap \mathcal{C}([0,T],%
\mathbb{R}
),$ $j=1,...,n$. The process $N:=(N^{1},...,N^{n})$ is defined by 
\begin{equation}
N_{t}^{j}=b_{t}^{j}+\displaystyle\displaystyle\int_{0}^{t}\sigma
_{s}^{j}\delta X_{s}^{j},\text{ }t\in \lbrack 0,T], j=1,...,n.  \label{N}
\end{equation}

Let $t_{0}\geqslant 0$ be fixed, and denote by $\mathbb{L}^{2}(\mathbb{F},%
\mathbb{%
\mathbb{R}
}^{n})$ the set of $\mathbb{F}$-adapted $\mathbb{%
\mathbb{R}
}^{n}$-valued processes $Z$ such that $\mathbb{E}(\int_{t_{0}}^{T}\mid
Z_{t}\mid ^{2}dt)<\infty .$ We consider the linear BDSE for the processes $%
Y=(Y_{t},$ $t\in \lbrack t_{0},T])\in \mathbb{L}^{2}(\mathbb{F},\mathbb{%
\mathbb{R}
})$ and $Z=((Z_{t}^{1},...,Z_{t}^{n}),$ $t\in \lbrack t_{0},T])\in \mathbb{L}%
^{2}(\mathbb{F},\mathbb{%
\mathbb{R}
}^{n})$ given by

\begin{equation}
Y_{t}=g(N_{T})-\int_{t}^{T}[f(s)+A_{1}(s)Y_{s}-A_{2}(s)Z_{s}]ds+%
\int_{t}^{T}Z_{s}\delta X_{s},\text{ \ }t\in \lbrack t_{0},T],  \label{BSDE}
\end{equation}

where the real-valued functions $g,f,A_{1}$ and the $%
\mathbb{R}
^{n}-$valued function $A_{2}$ are supposed to be known and the integral $%
\int_{t}^{T}Z_{s}\delta X_{s}$ is defined as a divergence integral and will
be studied in Sections 2 and 3. (\ref{BSDE}) is associated to the following
second order linear PDE with terminal condition%
\begin{equation}
\frac{\partial u}{\partial t}(t,x)=-\frac{1}{2}\sum_{j=1}^{n}\frac{d
}{dt}Var( N_{t}^{j})\frac{\partial ^{2}u}{\partial x_{j}^{2}%
}(t,x)+\sum_{j=1}^{n}[\sigma _{t}^{j}A_{2}^{j}(t)-\frac{d}{dt}b^{j}_t]%
\frac{\partial u}{\partial x_{j}}(t,x)+A_{1}(t)u(t,x)+f(t),  \label{PDE}
\end{equation}%
\begin{equation*}
u(T,x)=g(x),\ (t,x)\in \lbrack t_{0},T)\times \mathbb{R}^{n}.\text{ \ \ \ \
\ \ \ \ \ \ \ \ \ \ \ \ \ \ \ \ \ \ \ \ \ \ \ \ \ \ \ \ \ \ \ \ \ \ \ \ \ \
\ \ \ \ \ \ \ \ \ \ \ \ \ \ \ \ \ \ \ \ \ \ \ }
\end{equation*}

By means of the It\^{o} formula of Section 3 we show that 
\begin{equation}
Y_{t}:=u(t,N_{t})\ \text{and}\ Z_{t}^{j}:=-\sigma _{t}^{j}\frac{\partial u}{%
\partial x_{j}}(t,N_{t}),\ j=1,...,n  \label{YetZ}
\end{equation}%
is a solution of (\ref{BSDE})$.$ (\ref{PDE}) is solved explicitely in
Section 4.

A possible application of this result is to formulate and to solve, in the
context of Volterra processes, the well known problem of self-financing
trading strategies against a positive contingence claim $\xi $ at a time $T$%
. Volterra processes are in fact more flexible than brownian motion to
modelize the different kinds of uncertainties that occur on financial
markets. Let $P_{t}^{0}$ be the price at time $t\in \lbrack t_{0},T]$ of a
riskless asset that is governed by the equation $%
dP_{t}^{0}=P_{t}^{0}r_{t}dt, $ where $r_{t}$ is the short rate, and let $%
P_{t}^{j},$ $j=1,...,n$ be the prices of $n$ risky assets that are governed
by the equations $dP_{t}^{j}=P_{t}^{j}(b_{t}^{j}dt+\sigma _{t}^{j}\delta
X_{t}^{j}),$ where $b_{t}^{j}$ is the appreciation rate and $\sigma _{t}^{j}$
is the volatility. Then $\theta _{t}^{j}:=$ ($b_{t}^{j}-r_{t})/\sigma
_{t}^{j}$ is called the risk premium for the risky asset $j.$ The problem is
to find the part $\pi _{t}^{j}$ of the wealth $V_{t}$ to invest at time $t$
in the asset $j,$ $j=1,...,n.$ Of course this decision can only be based on
the current information $\mathcal{F}_{t},$ i.e. the processes $\pi _{t}^{j}$
and $\pi _{t}^{0}:=V_{t}-\sum_{j=1}^{n}\pi _{t}^{j}$ must be $\mathcal{F}_{t}
$-adapted. In the model of Harrison and Pliska (\cite{HP}) $V_{t}$ satisfies
the equation $dV_{t}=r_{t}V_{t}dt+\sum_{j=1}^{n}(b_{t}^{j}-r_{t})\pi
_{t}^{j}dt+\sum_{j=1}^{n}\pi _{t}^{j}\sigma _{t}^{j}\delta X_{t}^{j}$ with
final condition $V_{T}=\xi .$ If $\xi $ is of the form $g(N_{T}),$ ($%
V_{t},\pi _{t}^{1},...,\pi _{t}^{n}),$ $t\in \lbrack t_{0},T],$ is the
solution of (\ref{BSDE}) with $A_{1}(s)=r_{s},$ and $A_{2}^{j}(s)=\theta
_{s}^{j}\sigma _{s}^{j},$ $j=1,...,n.$

\section{On the divergence integral for Gaussian Volterra processes}

%
The kernel $K$ in (\ref{X}) defines a linear operator in $L^{2}([0,T])$
given by $(K\sigma )_{t}=\int_{0}^{t}K(t,s)\sigma _{s}ds,$ $\sigma \in
L^{2}([0,T]).$ Let $\mathcal{E}$ be the set of step functions of $[0,T]$,
and let $K_{T}^{\ast }:\mathcal{E}\rightarrow L^{2}([0,T])$ be defined by 
\begin{equation*}
(K_{T}^{\ast }\sigma )_{u}:=\int_{u}^{T}\sigma _s\frac{\partial K}{\partial
s}(s,u)ds.
\end{equation*}%
The operator $K_{T}^{\ast }$ is the adjoint of $K$ (\cite{AMN}, Lemma 1).%
\newline
%

\textbf{Remarks:} a) For $s>t,$ we have $(K_{T}^{\ast }\sigma
1_{[0,t]})_{s}=0,$ and we will denote $(K_{T}^{\ast }\sigma 1_{[0,t]})_{s}$
by $(K_{t}^{\ast }\sigma )_{s}$ where $K_{t}^{\ast }$ is the adjoint of the
operator $K$ in the interval $[0,t].$\newline

b) If $K(u,u)=0$ for all $u\in \lbrack 0,T],$ $(K_{T}^{\ast
}1_{[0,r]})_{u}=K(r,u)$ for $u<r.$ Indeed, if $u\leqslant r,$ we have\newline
\begin{equation*}
(K_{T}^{\ast }1_{[0,r]})_{u}=\displaystyle\int_{u}^{T}1_{[0,r]}(s)\frac{%
\partial K}{\partial s}(s,u)ds=\displaystyle\int_{u}^{r}\frac{\partial K}{%
\partial s}(s,u)ds.
\end{equation*}

Therefore %
%
%
%
\begin{align*}
R(t,s)& =\mathbb{E}\Big[X_{t}X_{s}\Big]=\displaystyle%
\int_{0}^{inf(t,s)}(K_{T}^{\ast }1_{[0,t]})_{u}(K_{T}^{\ast }1_{[0,s]})_{u}du
\\
& =<K_{T}^{\ast }1_{[0,t]},K_{T}^{\ast }1_{[0,s]}>_{L^{2}([0,T])}.
\end{align*}

For $\sigma ,\widetilde{\sigma }\in \mathcal{E}$ this may be extended to%
\begin{equation*}
X(\sigma ):=\displaystyle\int_{0}^{t}(K_{t}^{\ast }\sigma )_{s}dW_{s}\text{
\ and \ }E\Big[X(\sigma )X(\widetilde{\sigma })\Big]=<K_{T}^{\ast }\sigma
,K_{T}^{\ast }\widetilde{\sigma }>_{L^{2}([0,T])}.
\end{equation*}

Let $\mathcal{H}$ be the closure of the linear span of the indicator
functions $1_{[0,t]},$ $t\in \lbrack 0,T]$ with respect to the scalar product%
\begin{equation*}
<1_{[0,t]},1_{[0,s]}>_{\mathcal{H}}:=<K_{T}^{\ast }1_{[0,t]},K_{T}^{\ast
}1_{[0,s]}>_{L^{2}([0,T])}.
\end{equation*}%
The operator $K_{T}^{\ast }$ is an isometry between $\mathcal{H}$ and a
closed subspace of $L^{2}([0,T]),$ and $\parallel \cdot \parallel _{\mathcal{%
H}}$ is a semi-norm on $\mathcal{H}.$ Furthermore, for $\varphi ,\psi \in 
\mathcal{H}$, 
%
\begin{equation*}
<K_{T}^{\ast }\varphi ,K_{T}^{\ast }\psi >_{L^{2}([0,T])}=\displaystyle%
\int_{0}^{T}\displaystyle\int_{0}^{T}\Big(\displaystyle\int_{0}^{inf(r,s)}%
\frac{\partial K}{\partial r}(r,t)\frac{\partial K}{\partial s}(s,t)dt\Big)%
\varphi _{r}\psi _{s}dsdr.
\end{equation*}

For further use let

\begin{eqnarray*}
\phi (r,s) &:&=\displaystyle\int_{0}^{inf(r,s)}\frac{\partial K}{\partial r}%
(r,t)\frac{\partial K}{\partial s}(s,t)dt,\text{ }r\neq s,\text{\ } \\
\widetilde{\phi }(r,s) &:&=\displaystyle\int_{0}^{inf(r,s)}\Big | \frac{%
\partial K}{\partial r}(r,t)\Big | \Big |\frac{\partial K}{\partial s}(s,t)\Big |
dt,\text{ }r\neq s.\text{\ }
\end{eqnarray*}%
Note that $\phi (r,s)=\partial ^{2}/\partial s\partial rR(r,s)$ ($r\neq s)$ (%
$\phi $ may be infinite on the diagonal $r=s).$ 
%
%
%
%
%
%
%
Let $\mid \mathcal{H}\mid $ be the closure of the linear span of indicator
functions with respect to the semi-norm given by

\begin{align*}
\parallel \varphi \parallel _{\mid \mathcal{H}\mid }^{2}&=\displaystyle%
\int_{0}^{T}\left( \int_{t}^{T}\mid \varphi _{r}\mid \Big | \frac{\partial K}{%
\partial r}(r,t)\Big | dr\right) ^{2}dt \\
&=2\displaystyle\int_{0}^{T}dr\displaystyle\int_{0}^{r}ds\widetilde{\phi }%
(r,s)\mid \varphi _{r}\mid \mid \varphi _{s}\mid .
\end{align*}

We briefly recall some basic elements of the stochastic calculus of
variations with respect to $X.$ We refer to \cite{Nualart} for a more
complete presentation. Let $\mathcal{S}$ be the set of random variables of
the form $F=f(X(\varphi _{1}),....,X(\varphi _{n}))$, where $n\geqslant 1$, $%
f\in \mathcal{C}_{b}^{\infty }(\mathbb{R}^{n})$ ($f$ and its derivatives are
bounded) and $\varphi _{1},...,\varphi _{n}\in \mathcal{H}$. The derivative
of $F$ 
\begin{equation*}
\begin{aligned} D^{X}F &:= \sum_{j=1}^{n}\frac{\partial f}{\partial
x_{j}}(X(\varphi_1),....,X(\varphi_n))\varphi_{j}, \end{aligned}
\end{equation*}%
\quad is an $\mathcal{H}$-valued random variable, and $D^{X}$ is a closable
operator from $L^{p}(\Omega )$ to $L^{p}(\Omega ;\mathcal{H})$ for all $%
p\geqslant 1$. We denote by $\mathbb{D}_{1,p}^{X}$ the closure of $\mathcal{S%
}$ with respect to the semi-norm 
\begin{equation}
\begin{aligned}
\|F\|_{1,p}^{p}&=\mathbb{E}|F|^{p}+\mathbb{E}\|D^{X}F\|_{\mathcal{H}}^{p}.
\end{aligned}
\label{normeD1p}
\end{equation}%
\newline
We denote by $Dom (\delta^X)$ the subset of $L^{2}(\Omega, \mathcal{H})$ composed of those elements $u$ for which there exists a positive constant $c$ such that
\begin{equation}
\Big | \mathbb{E} \Big[ <D^XF, u>_{\mathcal{H}}\Big]\Big | \leqslant c\sqrt{\mathbb{E}[F^2]},\, \text{for all $F\in\mathbb{D}_{1,2}^{X}.$}
\label{Domaiendelta}
\end{equation}
 For $u\in L^{2}(\Omega ;\mathcal{H}%
)$ in $Dom (\delta^X),$ $\delta ^{X}(u)$ is the element in $%
L^{2}(\Omega )$ defined by the duality relationship 
\begin{equation}
\mathbb{E}\Big[F\delta ^{X}(u)\Big]=\mathbb{E}\Big[<D_{\cdot
}^{X}F, u_{\cdot }>_{\mathcal{H}}\Big],\text{ }F\in \mathbb{D}_{1,2}^{X}.
\label{dualité}
\end{equation}

We also use the notation $\int_{0}^{T}u_{t}\delta X_{t}$ for\ $\delta
^{X}(u).$ A\ class of processes that belong to the domain of $\delta ^{X}$
is given as follows: let $\mathcal{S}^{\mathcal{H}}$ be the class of $%
\mathcal{H}$-valued random variables $u=\sum_{j=1}^{n}F_{j}h_{j}$ ($F_{j}\in 
\mathcal{S},$ $h_{j}\in \mathcal{H}).$ \\
In the same way $\mathbb{D}_{1,p}^{X}(\mid \mathcal{H\mid })$ is defined as
the completion of $\mathcal{S}^{\mid \mathcal{H\mid }}$ under the semi-norm \\  
$$\parallel u\parallel _{1,p,\mid \mathcal{H\mid }}^{p}:=\mathbb{E}\parallel u\parallel_{\mid\mathcal{H}\mid}^p+\mathbb{E}\parallel D^Xu\parallel_{\mid\mathcal{H}\mid\otimes\mid\mathcal{H}\mid}^p,$$ where%
\begin{equation}
\parallel D^{X}u\parallel _{\mid \mathcal{H\mid \otimes \mid H\mid }%
}^{2}=\int_{[0,T]^{4}}\mid D_{s}^{X}u_{t}\mid \mid D_{t^{\prime
}}^{X}u_{s^{\prime }}\mid \widetilde{\phi }(s, s^{\prime })\widetilde{\phi }%
(t, t^{\prime })dsdtds^{\prime }dt^{\prime }.
\label{normeDXu}
\end{equation}
The space $\mathbb{D}_{1, 2}^{X}(\mid\mathcal{H}\mid)$ is included in the domain of $\delta ^X$, and we have, for $u\in\mathbb{D}_{1, 2}^{X}(\mid\mathcal{H}\mid),$
$$\mathbb{E}\Big(\delta^X(u)^2\Big)\leqslant \mathbb{E}\parallel u\parallel_{\mid\mathcal{H}\mid}^2+\mathbb{E}\parallel D^X u\parallel_{\mid\mathcal{H}\mid\otimes\mid\mathcal{H}\mid}^2.$$
For $F\in \mathcal{S},$ let%
\begin{equation}
\mathbb{D}_{s}^{X}F:=\displaystyle\int_{0}^{T}\phi (s, t)D_{t}^{X}Fdt,\,
s\in \lbrack 0,T].
\label{Ddoublebar}
\end{equation}
Then (\ref{Ddoublebar}) implies
\begin{equation}
\int_{[0,T]^{2}}\mathbb{D}_{s}^{X}u_{t}\mathbb{D}_{t}^{X}v_{s}dsdt=\int_{[0,T]^{4}}D_{s}^{X}u_{t}D_{t^{\prime }}^{X}v_{s^{\prime }}\phi
(s, s^{\prime })\phi (t, t^{\prime })dsdtds^{\prime }dt^{\prime }.
\label{DXuDXv}
\end{equation}

\begin{prop}
Let $f$,$g\in $ $\mathbb{D}_{1,2}^{X}(\mid \mathcal{H\mid }).$ Then the
integrals $\delta ^{X}(f)$ and $\delta ^{X}(g)$ exist in $L^{2}(\Omega )$
and 
\begin{equation}
\mathbb{E}\Big [\delta ^{X}(f)\delta ^{X}(g)\Big]=\mathbb{E}<f,g>_{\mathcal{H%
}}+\displaystyle\int_{0}^{T}ds\displaystyle\int_{0}^{T}dt\mathbb{E}\Big[%
\mathbb{D}_{t}^{X}f_{s}\mathbb{D}_{s}^{X}g_{t}\Big].  \label{cov}
\end{equation}
\label{propcovariance}
\end{prop}

\textbf{Remark. }With the choice $f=g,$ Proposition \ref{propcovariance}
implies $\mathbb{D}_{1,2}^{X}(\mid \mathcal{H\mid })\subset dom$($\delta
^{X}).$ In fact, for $f\in $ $\mathbb{D}_{1,2}^{X}(\mid \mathcal{H\mid }),$%
\begin{equation*}
\mathbb{\parallel }f\mathbb{\parallel }_{\mathcal{H}}\leqslant \mathbb{%
\parallel }f\mathbb{\parallel }_{\mid \mathcal{H\mid }}\text{ and }%
\int_{[0,T]^{2}}\mathbb{D}_{s}^{X}f_{t}\mathbb{D}_{t}^{X}f_{s}dsdt\leqslant%
\parallel D^{X}f\parallel _{\mid \mathcal{H\mid \otimes \mid H\mid }}^{2}.
\end{equation*}

For the proof we need the following lemma.

\begin{lem}
Let $G\in \mathbb{D}_{1,2}^{X}$ and $f$ $\in dom(\delta ^{X}).$ Then 
\begin{equation*}
\mathbb{E}\Big [\delta ^{X}(f)G\Big]=\displaystyle\int_{0}^{T}\mathbb{E}\Big[%
f_{t}\mathbb{D}_{t}^{X}G\Big]dt.
\end{equation*}
\end{lem}

\quad Proof of the Lemma. 
\begin{align*}
\mathbb{E}\Big [\delta ^{X}(f)G\Big]& =\mathbb{E}\Big[<D_{\cdot
}^{X}G,f_{\cdot }>_{\mathcal{H}}\Big]=\mathbb{E}\Big[<(K_{T}^{\ast
}D^{X}G)_{\cdot },(K^{\ast }f)_{\cdot }>_{L^{2}([0,T])}\Big] \\
& =\mathbb{E}\Big[\displaystyle\int_{0}^{T}dt\displaystyle%
\int_{t}^{T}du_{1}f_{u_{1}}\frac{\partial K}{\partial u_{1}}(u_{1},t)%
\displaystyle\int_{t}^{T}du_{2}D_{u_{2}}^{X}G\frac{\partial K }{\partial u_{2}}%
(u_{2},t)\Big] \\
& =\mathbb{E}\Big[\displaystyle\int_{0}^{T}du_{1}f_{u_{1}}\displaystyle%
\int_{0}^{T}du_{2}D_{u_{2}}^{X}G\phi (u_{1},u_{2})\Big]=\mathbb{E}\Big[%
\displaystyle\int_{0}^{T}f_{u_{1}}\mathbb{D}_{u_{1}}^{X}Gdu_{1}\Big]%
.\ \ \ \ \blacksquare
\end{align*}

Proof of Proposition \ref{propcovariance}. Let first $f,g\in \mathcal{S}%
^{\mid \mathcal{H\mid }}.$ If $f=\sum_{j=1}^{n}F_{j}h_{j},$ $F_{j}\in 
\mathcal{S},$ $h_{j}\in \mid \mathcal{H\mid },$ lemma 2.2 in \cite{Nualart}
shows that 
\begin{equation*}
\delta ^{X}(f)=\sum_{j=1}^{n}\left( F_{j}X(h_{j})-\langle
D^{X}F_{j},h_{j}\rangle _{\mathcal{H}}\right) \in \mathcal{S}.
\end{equation*}

A straightforward calculation shows $D_{.}^{X}\delta ^{X}(f)\in \mathcal{S}%
^{\mid \mathcal{H\mid }}$ and $\delta ^{X}(D_{.}^{X}f)\in \mathcal{S}^{\mid 
\mathcal{H\mid }}.$ Moreover, by means of the notion of directional
derivatives (\cite{Nualart}, p.27), one shows%
\begin{equation*}
\left\langle f_{\cdot },D_{\cdot }^{X}(\delta ^{X}(g))\right\rangle _{%
\mathcal{H}}=\left\langle f_{\cdot },\delta ^{X}(D_{\cdot }^{X}g)+g_{\cdot
}\right\rangle _{\mathcal{H}}.
\end{equation*}

Therefore

\begin{equation*}
\mathbb{E}\Big [\delta ^{X}(f)\delta ^{X}(g)\Big]=\mathbb{E}\left\langle
f_{\cdot },D_{\cdot }^{X}(\delta ^{X}(g))\right\rangle _{\mathcal{H}}=%
\mathbb{E}<f,g>_{\mathcal{H}}+\mathbb{E}\left\langle f_{\cdot },\delta
^{X}(D_{\cdot }^{X}g)\right\rangle _{\mathcal{H}},
\end{equation*}
Moreover 
\begin{eqnarray*}
\mathbb{E}\left\langle f_{\cdot },\delta ^{X}(D_{\cdot }^{X}g)\right\rangle
_{\mathcal{H}} &=&\mathbb{E}\left[ \int_{0}^{T}(K_{T}^{\ast }f)_{u}\left(
K_{T}^{\ast }\delta ^{X}(D_{\cdot }^{X}g)\right) _{u}du\right]  \notag \\
&=&\mathbb{E}\left[ \int_{0}^{T}du\int_{u}^{T}dsf_{s}\frac{\partial K}{%
\partial s}(s,u)\int_{u}^{T}dr\delta ^{X}(D_{r}^{X}g)\frac{\partial K}{%
\partial r}(r,u)\right]  \notag \\
&=&\mathbb{E}\left[ \int_{0}^{T}ds\int_{0}^{T}drf_{s}\delta
^{X}(D_{r}^{X}g)\left( \int_{0}^{r\wedge s}\frac{\partial K}{\partial s}(s,u)%
\frac{\partial K}{\partial r}(r,u)du\right) \right]  \notag \\
&=&\mathbb{E}\left[ \int_{0}^{T}ds\int_{0}^{T}drf_{s}\delta
^{X}(D_{r}^{X}g)\phi (s,r)\right] =\int_{0}^{T}ds\int_{0}^{T}dr\mathbb{E}%
[f_{s}\delta ^{X}(D_{r}^{X}g)]\phi (s,r)  \notag \\
&=&\int_{0}^{T}ds\int_{0}^{T}dr\mathbb{E}\Big[%
<D_{r}^{X}g_{.},D_{.}^{X}f_{s}>_{\mathcal{H}}\Big]\phi (s,r)  \notag \\
&=&\int_{0}^{T}ds\int_{0}^{T}dr\mathbb{E}\Big[<K_{T}^{\ast
}D_{r}^{X}g_{.},K_{T}^{\ast }D_{.}^{X}f_{s}>_{L^{2}([0,T])}\Big]\phi (s,r) \\
&=&\int_{0}^{T}ds\mathbb{E}\Big[<K_{T}^{\ast }\int_{0}^{T}D_{r}^{X}g_{.}\phi
(s,r)dr,K_{T}^{\ast }D_{.}^{X}f_{s}>_{L^{2}([0,T])}\Big] \\
&=&\int_{0}^{T}ds\mathbb{E}<\int_{0}^{T}D_{r}^{X}g_{.}\phi
(s,r)dr,D_{.}^{X}f_{s}>_{\mathcal{H}}=\int_{0}^{T}ds\mathbb{E}\left[
f_{s}\delta ^{X}\left( \int_{0}^{T}D_{r}^{X}g_{.}\phi (s,r)dr\right) \right] 
\notag \\
&=&\int_{0}^{T}ds\mathbb{E}\left[ f_{s}\delta ^{X}\left( \mathbb{D}%
_{s}^{X}g_{.}\right) \right] =\int_{0}^{T}ds\int_{0}^{T}dt\mathbb{E}\left[ 
\mathbb{D}_{s}^{X}g_{t}\mathbb{D}_{t}^{X}f_{s}\right] ,  \notag
\end{eqnarray*}

where we have used the preceding lemma in the last equality. This proves (%
\ref{cov}) for $f,g\in \mathcal{S}^{\mid \mathcal{H\mid }}.$

Let now $f\in $ $\mathbb{D}_{1,2}^{X}(\mid \mathcal{H\mid }).$ There exists
a sequence $(f_{n})_{n\in 
\mathbb{N}
}\subset \mathcal{S}^{\mid \mathcal{H\mid }}$ such that, as $n\rightarrow
\infty ,$ $f_{n}\rightarrow f$ in $L^{2}(\Omega ,\mathcal{H})$ and $%
D^{X}f_{n}\rightarrow D^{X}f$ in $L^{2}(\Omega ,\mathcal{H\otimes H}).$
Therefore $\delta ^X(f_{n})$ converges in $L^{2}(\Omega ).$ By proposition
1.3.6 in \cite{Nualart} the limit is $\delta^X (f).$\ \ \ \ $\blacksquare $%
\quad\

\section{It\^{o} formula}

Let $F\in \mathcal{C}^{1,2}([0,T]\times 
\mathbb{R}
^{n})$ and suppose that 
\begin{equation}
\begin{aligned} \max \Big( \Big | F(t,x)\Big |,\Big |\frac{\partial
F}{\partial t}(t,x)\Big |,\Big |{\frac{\partial F}{\partial x_j}}(t,x)\Big
|,\Big |{\frac{\partial^{2}F}{\partial x_j^{2}}(t,x)}\Big |, j=1,...,n \Big)
&\leqslant ce^{\lambda \mid x\mid^{2}} \end{aligned}  \label{F}
\end{equation}%
for all $t\in \lbrack 0,T]$ and $x\in \mathbb{R}^{n},$ where $c,\lambda $
are positive constants such that $\lambda <\frac{1}{4}\underset{i}{\min }(%
\underset{t\in \lbrack 0,T]}{\sup }Var(N_{t}^{i}))^{-1}.$ This
implies 
\begin{eqnarray*}
\mathbb{E}\Big |F(t, N_{t})\Big |^{2} &\leqslant &c^{2}\mathbb{E}exp(2\lambda
|N_{t}|^{2}) \\
&=&\frac{c^{2}}{(2\pi )^{n/2}}\prod\limits_{j=1}^{n}\frac{1}{(Var%
(N_{t}^{j}))^{1/2}}\displaystyle\int exp\left( 2\lambda x_{j}^{2}-\frac{1%
}{2}\,\frac{(x_{j}-b_t^j)^{2}}{Var(N_{t}^{j})}\right) dx_{j}<\infty ,
\end{eqnarray*}

%
%
%
%
%
%
%
%
%
%
%
%
%
%
%
%
%
%
and the same property holds for $\partial /\partial tF(t, x),$ $\partial
/\partial x_{j}F(t, x)$ and $\partial ^{2}/\partial x_{j}^{2}F(t, x),$ $%
j=1,..., n$.

\begin{thm}
Let $N\ $be given by (\ref{N}), and suppose that, for $j=1,...,n,$ the
kernels $K^{j}$ of $X^{j}$ satisfy (H1) and (H2), $b^{j}\in \mathcal{C}%
^{1}((0,T),%
\mathbb{R}
)\cap \mathcal{C}([0,T],%
\mathbb{R}
)$, and $\sigma =\{\sigma _{t}^{j},$ $t\in \lbrack 0,T],$ $j=1,...,n\}$ is
bounded. If $F\in $ $\mathcal{C}^{1,2}([0,T]\times \mathbb{R}^{n})$
satisfies $(\ref{F})$, $\partial /\partial x_{j}F(\cdot ,N_{\cdot })\in 
\mathbb{D}_{1,2}^{X^{j}}(\mid \mathcal{H}^{j}\mathcal{\mid }),$ $j=1,...,n$
and, for all $t\in \lbrack 0,T],$%
\begin{align}
F(t,N_{t})& =F(0,0)+\int_{0}^{t}\frac{\partial F}{\partial s}%
(s,N_{s})ds+\sum_{j=1}^{n}\int_{0}^{t}\frac{\partial F}{\partial x_{j}}%
(s,N_{s})\Big(\frac{d}{ds}b^{j}_s ds+\sigma _{s}^{j}\delta X_{s}^{j}\Big) 
\notag \\
& +\frac{1}{2}\sum_{j=1}^{n}\int_{0}^{t}\frac{\partial ^{2}F}{\partial
x_{j}^{2}}(s,N_{s})\frac{d}{ds}Var(N_{s}^{j})ds.
\label{Itodiagonal}
\end{align}%
\label{FormulaItodiagonal}
\end{thm}

\textbf{Remarks. a) }A more general model for $N$ than in Section 1.2 is

\begin{equation}
\widetilde{N}_{t}^{i}=b_{t}^{i}+\sum_{j=1}^{n}\int_{0}^{t}\widetilde{\sigma }%
_{s}^{i,j}\delta X_{s}^{j},\text{ }i=1,...,n,
\end{equation}%
where $\widetilde{\sigma }=(\widetilde{\sigma }_{t}^{i,j},$ $i,j=1,...,n)$
is a matrix of bounded functions $\widetilde{\sigma }^{i,j}$ defined on $%
[0,T].$ Let $\widetilde{N}=(\widetilde{N}^{1},...,\widetilde{N}^{n}).$ The
components of $\widetilde{N}$ are dependent since $\widetilde{N}^{i}$
depends not only on a single random perturbation $X^{i}$, but on the others $%
X^{j}$ ($j\neq i$) as well. The model of Section 1.2\ may be recovered by
choosing the matrix $\widetilde{\sigma }$ diagonal with\ functions $\sigma
^{j}:=\widetilde{\sigma }^{j,j}$ in the diagonal.\emph{\ }An It\^{o} formula
can be shown for $F(t,\widetilde{N_{t}})$ too, but, instead of the variances
of $N,$ the covariances of $\widetilde{N}$ appear now in the second order
term. It reads 
\begin{align}
F(t,\widetilde{N}_{t})& =F(0,0)+\int_{0}^{t}\frac{\partial F}{\partial s}(s,%
\widetilde{N}_{s})ds+\sum_{i=1}^{n}\int_{0}^{t}\frac{\partial F}{\partial
x_{i}}(s,\widetilde{N}_{s})\Big(\frac{d}{ds}b^{i}_s ds+\sum_{j=1}^{n}%
\widetilde{\sigma }_{s}^{i,j}\delta X_{s}^{j}\Big)  \notag \\
& +\frac{1}{2}\sum_{i=1}^{n}\sum_{j=1}^{n}\int_{0}^{t}\frac{\partial ^{2}F}{%
\partial x_{i}\partial x_{j}}(s,\widetilde{N}_{s})\frac{d}{ds}Cov(%
\widetilde{N}_{s}^{i},\widetilde{N}_{s}^{j})ds.  \label{Itondim}
\end{align}%
Two problems arise with this model. The first is to find a positive constant 
$\lambda $ in the growth condition (\ref{F}) that implies $\mathbb{E}\Big |%
F(t,\widetilde{N_{t}})\Big |^{2}<\infty .$ It seems therefore that (\ref%
{Itondim}) can be shown for $F\in $ $\mathcal{C}_{b}^{1,2}([0,T]\times 
\mathbb{R}^{n})$ only, i.e. if $F$ and its derivatives are bounded, a
hypothesis that is in general too restrictive for an application to the
solution $u$ of the PDE (\ref{PDE}). The second problem is that, in order to
get the solution of the BSDE (\ref{BSDE}) with $\widetilde{N},$ the
derivatives of the variances of $N$ in the second order term of (\ref{PDE})
should be replaced by the derivatives of the covariances of $\widetilde{N}$,
and it would be necessary to assume that the matrix ($\frac{d}{dt}Cov(%
\widetilde{N}_{t}^{i},\widetilde{N}_{t}^{j}),$ $i,j=1,...,n)$ is positive
definite, a hypothesis which seems to be restrictive (see also the remark
below). For this reason we prefer to consider BSDEs with $N$ given by (\ref%
{N}). We notice that 
\begin{equation*}
Var(N_{t}^{j})=\mathbb{E}\left( \int_{0}^{t}\left( K_{t}^{\ast
,j}\sigma ^{j}\right) _{s}\delta W_{s}^{j}\right) ^{2}=\int_{0}^{t}\Big(%
K_{t}^{\ast ,j}\sigma ^{j}\Big)_{s}^{2}ds,\text{ and }\frac{d}{dt}Var%
(N_{t}^{j})=2\sigma _{t}^{j}\mathbb{D}_{t}^{j}N_{t}^{j},\text{ \ \ }
\end{equation*}

\textbf{b) }In the model where $\widetilde{\sigma }$ is diagonal, $X^{i}$
and $X^{j}$ are defined with independent brownian motions $W^{i}$ and $W^{j}$
if $i\neq j.$ This differs from the model where all the Volterra processes $%
\overline{X}^{j}$ are defined with the same brownian motion $\overline{W}$
as follows:%
\begin{equation*}
\overline{X}_{t}^{j}=\int_{0}^{t}K^{j}(t,s)\delta \overline{W}_{s},\text{ }%
\overline{N}_{t}^{j}=b_{t}^{j}+\int_{0}^{t}\sigma _{s}^{j}\delta \overline{X}%
_{s}^{j},\text{ }t\in \lbrack 0,T],\text{ }j=1,...,n.
\end{equation*}

In this case the processes $\overline{N}^{j}$ are correlated, and the matrix
($\frac{d}{dt}Cov(\overline{N}_{t}^{i},\overline{N}_{t}^{j}),$ $%
i,j=1,...,n)$ is not diagonal. For $n=2$ a\ direct verification shows that
this matrix is negative (semi-)definite.\newline

Proof of Theorem \ref{FormulaItodiagonal}\newline
1. Let us first show that $\partial /\partial x_{j}F(\cdot ,N_{\cdot })\in
\delta ^{X^{j}}$ for all $j=1,...,n.$ For this we show that $\partial
/\partial x_{j}F(\cdot ,N_{\cdot })\in \mathbb{D}_{1,2}^{X^{j}}(\mid 
\mathcal{H}^{j}\mid ),$ where $\mathcal{H}^{j}$ is the space defined in
Section 2 with $X$ replaced by $X^{j}$. The terms $\phi ^{j}$ and $%
\widetilde{\phi ^{j}}$ refer now to the kernel $K^{j}$ of $X^{j}.$ The
constants in the inequalities below may vary from line to line. 
\begin{align*}
\mathbb{E}\Big\|\frac{\partial F}{\partial x_{j}}(.,N_{.})\Big\|_{\mathcal{H}%
^{j}}^{2}& =\mathbb{E}\displaystyle\int_{0}^{T}\Big(\int_{s}^{T}\Big|\frac{%
\partial F}{\partial x_{j}}(t,N_{t})\Big|\Big|\frac{\partial K^{j}}{\partial
t}(t,s)\Big|dt\Big)^{2}ds \\
& \leqslant c\mathbb{E}\displaystyle\int_{0}^{T}\Big(\int_{s}^{T}\exp
(\lambda \mid N_{t}\mid ^{2})(t-s)^{\alpha -1}\Big(\frac{t}{s}\Big)^{\beta
}dt\Big)^{2}ds
\end{align*}%
Applying H\"{o}lder's inequality to $1<p<\frac{1}{1-\alpha }<2$ and $q$ its
conjugate$,$ we get 
\begin{equation*}
\Big(\int_{s}^{T}\exp (\lambda \mid N_{t}\mid ^{2})(t-s)^{\alpha -1}\Big(%
\frac{t}{s}\Big)^{\beta }dt\Big)^{2}\leqslant \Big(\int_{s}^{T}\exp
(q\lambda \mid N_{t}\mid ^{2})dt\Big)^{\frac{2}{q}}\Big(%
\int_{s}^{T}(t-s)^{p(\alpha -1)}\Big(\frac{t}{s}\Big)^{p\beta }dt\Big)^{%
\frac{2}{p}}.
\end{equation*}%
Then%
\begin{equation*}
\mathbb{E}\Big(\int_{0}^{T}\exp (q\lambda \mid N_{t}\mid ^{2})dt\Big)^{\frac{%
2}{q}}\leqslant \mathbb{E}\Big(\int_{0}^{T}\exp (q\lambda \underset{t\in
\lbrack 0,T]}{\sup }\mid N_{t}\mid ^{2})dt\Big)^{\frac{2}{q}}\leqslant
T^{2/q}\prod\limits_{i=1}^{n}\mathbb{E}\exp (2\lambda \underset{t\in \lbrack
0,T]}{\sup }(N_{t}^{i})^{2})dt<\infty
\end{equation*}%
for $\lambda <\frac{1}{4}\underset{i}{\min }(\underset{t\in \lbrack 0,T]}{%
\sup }Var(N_{t}^{i}))^{-1}.$ Moreover, 
\begin{equation*}
\int_{0}^{T}\Big(\int_{s}^{T}(t-s)^{p(\alpha -1)}\Big(\frac{t}{s}\Big)%
^{p\beta }dt\Big)^{\frac{2}{p}}ds\leqslant T^{2(\alpha -1)+\frac{2}{p}%
+1}\beta (1-2\beta ,2(\alpha -1)+\frac{2}{p}+1).
\end{equation*}%
It remains to show that $\mathbb{E}\Big\|D^{X^{j}}\frac{\partial F}{\partial x_{j}}%
(.,N_{.})\Big\|_{\mid \mathcal{H}^{j}\mathcal{\mid \otimes \mid H}^{j}%
\mathcal{\mid }}<\infty .$

\begin{align*}
&\mathbb{E}\Big\| D^{X^{j}}\frac{\partial F}{\partial x_{j}}(.,N_{.})\Big\| _{\mid 
\mathcal{H}^{j}\mathcal{\mid \otimes \mid H}^{j}\mathcal{\mid }}^{2} \\
&=\mathbb{E}\int_{[0,T]^{4}}\Big| D_{s}^{X^{j}}\frac{\partial F}{\partial
x_{j}}(t,N_{t})\Big| \Big| D_{t^{\prime }}^{X^{j}}\frac{\partial F}{\partial
x_{j}}(s^{\prime },N_{s^{\prime }})\Big| \widetilde{\phi ^{j}}(s,s^{\prime })%
\widetilde{\phi ^{j}}(t,t^{\prime })dsdtds^{\prime }dt^{\prime } \\
&=\mathbb{E}\int_{[0,T]^{4}}\Big|\frac{\partial ^{2}F}{\partial x_{j}^{2}}%
(t,N_{t})\sigma _{s}^{j}\Big| \Big| \frac{\partial ^{2}F}{\partial x_{j}^{2}}%
(s^{\prime },N_{s^{\prime }})\sigma _{t^{\prime }}^{j}\Big| \widetilde{\phi
^{j}}(s,s^{\prime })\widetilde{\phi ^{j}}(t,t^{\prime })dsdtds^{\prime
}dt^{\prime } \\
& \leqslant c\int_{0}^{T}dt\int_{0}^{T}ds^{\prime }\left[ \mathbb{E}\left( 
\frac{\partial ^{2}F}{\partial x_{j}^{2}}(t,N_{t})\right) ^{2}+\mathbb{E}%
\left( \frac{\partial ^{2}F}{\partial x_{j}^{2}}(s^{\prime },N_{s^{\prime
}})\right) ^{2}\right] \displaystyle\int_{0}^{T}dt^{\prime }\widetilde{\phi
^{j}}(t,t^{\prime })\int_{0}^{T}ds\widetilde{\phi ^{j}}(s^{\prime },s).
\end{align*}

By (\ref{F}), 
\begin{equation*}
\mathbb{E}\left( \frac{\partial ^{2}F}{\partial x_{j}^{2}}(t,N_{t})\right)
^{2}\leqslant c\mathbb{E}\exp (2\lambda \mid N_{t}\mid ^{2})\leqslant
c^{\prime }\prod\limits_{i=1}^{n}(1-4\lambda Var%
(N_{t}^{i}))^{-1/2},
\end{equation*}

and by the choice of $\lambda $ the right side stays bounded in $t\in
\lbrack 0,T].$ The finiteness of the remaining integrals follows from \emph{%
(H1), (H2)} applied to $K^{j}.$

2. We proceed now to the outline of the proof of the It\^{o} formula. Let 
\begin{equation*}
N_{t}^{j,\varepsilon }=b_{t}^{j}+\int_{0}^{t}\Big(\int_{s}^{t}\sigma _{r}^{j}%
\frac{\partial K^{j}}{\partial r}(r+\varepsilon ,s)dr\Big)\delta
W_{s}^{j}=b_{t}^{j}+\int_{0}^{t}\Big(\int_{0}^{r}\sigma _{r}^{j}\frac{%
\partial K^{j}}{\partial r}(r+\varepsilon ,s)\delta W_{s}^{j}\Big)dr,~~t\leqslant
T-\varepsilon .
\end{equation*}%
Then, $F(t,N_{t}^{1,\varepsilon },...,N_{t}^{n,\varepsilon })$ has locally
bounded variation, and we can write \newline
\begin{eqnarray*}
dF(t,N_{t}^{1,\varepsilon },...,N_{t}^{n,\varepsilon }) &=&\frac{\partial F}{%
\partial t}(t,N_{t}^{1,\varepsilon },...,N_{t}^{n,\varepsilon })dt+%
\displaystyle\sum_{i=1}^{n}\frac{\partial F}{\partial x_{i}}%
(t,N_{t}^{1,\varepsilon },...,N_{t}^{n,\varepsilon })dN_{t}^{i,\varepsilon }
\\
&=&\left( \frac{\partial F}{\partial t}(t,N_{t}^{\varepsilon
})+\sum_{j=1}^{n}\frac{\partial F}{\partial x_{j}}(t,N_{t}^{\varepsilon
})\frac{d}{dt}b^{j}_t\right) dt+\displaystyle\sum_{j=1}^{n}\frac{\partial F%
}{\partial x_{j}}(t,N_{t}^{\varepsilon })\displaystyle\sigma
_{t}^{j}\int_{0}^{t}\frac{\partial K^{j}}{\partial t}(t+\varepsilon
,s)\delta W_{s}^{j}dt
\end{eqnarray*}%
with the notation $N^{\varepsilon }=(N^{1,\varepsilon },...,N^{n,\varepsilon
}).$ Furthermore, 
\begin{align*}
\frac{\partial F}{\partial x_{j}}(t,N_{t}^{\varepsilon })\displaystyle\sigma
_{t}^{j}\int_{0}^{t}\frac{\partial K^{j}}{\partial t}(t+\varepsilon
,s)\delta W_{s}^{j}& =\displaystyle\sigma _{t}^{j}\Big[\int_{0}^{t}\frac{%
\partial F}{\partial x_{j}}(t,N_{t}^{\varepsilon })\frac{\partial K^{j}}{%
\partial t}(t+\varepsilon ,s)\delta W_{s}^{j} \\
& +\int_{0}^{t}D_{s}^{W^{j}}\Big(\frac{\partial F}{\partial x_{j}}%
(t,N_{t}^{\varepsilon })\Big)\frac{\partial K^{j}}{\partial t}(t+\varepsilon
,s)ds\Big],
\end{align*}

where 
\begin{equation*}
D_{s}^{W^{j}}\Big(\frac{\partial F}{\partial x_{j}}(t,N_{t}^{\varepsilon })%
\Big)=\frac{\partial ^{2}F}{\partial x_{j}^{2}}(t,N_{t}^{\varepsilon
})D_{s}^{W^{j}}N_{t}^{j,\varepsilon }=\frac{\partial ^{2}F}{\partial
x_{j}^{2}}(t,N_{t}^{\varepsilon })\int_{s}^{t}\sigma _{r}^{j}\frac{\partial
K^{j}}{\partial r}(r+\varepsilon ,s)dr.
\end{equation*}%
Therefore 
\begin{align*}
F(t,N_{t}^{\varepsilon })& =F(0,0)+\int_{0}^{t}\left( \frac{\partial F}{%
\partial s}(s,N_{s}^{\varepsilon })+\sum_{j=1}^{n}\frac{\partial F}{\partial
x_{j}}(s,N_{s}^{\varepsilon })\frac{d}{ds}b^{j}_s\right) ds \\
& +\sum_{j=1}^{n}\int_{0}^{t}\Big(\int_{s}^{t}\sigma _{r}^{j}\frac{\partial F%
}{\partial x_{j}}(r,N_{r}^{\varepsilon })\frac{\partial K^{j}}{\partial r}%
(r+\varepsilon ,s)dr\Big)\delta W_{s}^{j} \\
& +\frac{1}{2}\sum_{j=1}^{n}\int_{0}^{t}\frac{\partial ^{2}F}{\partial
x_{j}^{2}}(r,N_{r}^{\varepsilon })\frac{\partial }{\partial r}\Big(%
\int_{0}^{r}\Big(K_{r}^{\ast ,\varepsilon ,j}\sigma ^{j}\Big)_{s}^{2}ds\Big)%
dr.
\end{align*}%
We notice that the divergence integral in the second line above appears now
outside of the integral\ with respect to $dr.$ This permutation may be shown
with the definition by duality of the divergence integral. The integral
coincides, up to $\varepsilon $, with the divergence integral that appears
in the statement of the theorem. The last term coincides, up to $\varepsilon
,$ with the term in Remark a) after the theorem. It remains to show that the terms above converge in $%
L^{2}(\Omega )$ towards the terms in the statement of the theorem as $%
\varepsilon \rightarrow 0$. This can be done for each integral similarly as
in the proof of theorem 4 in \cite{AMN}. $\ \ \blacksquare $

\section{Solvability of linear BSDEs}

As mentioned in the introduction the aim is to apply the It\^{o} formula (%
\ref{Itodiagonal}) to the classical solution$\ u$\ of the PDE $(\ref{PDE})$
and to show that $Y$ and $Z$ defined by $(\ref{YetZ})$ satisfy the BSDE $(%
\ref{BSDE}).$ We will show later in this section that $u$\ in fact satisfies
the growth condition $(\ref{F})$ under a suitable growth condition on the
final condition in $(\ref{BSDE}).$\newline
The It\^{o} formula $(\ref{Itodiagonal}),$ applied to $u$ reads 
\begin{align}
u(t,N_{t})& =u(T,N_{T})-\int_{t}^{T}\frac{\partial u}{\partial s}%
(s,N_{s})ds-\sum_{j=1}^{n}\int_{t}^{T}\frac{\partial u}{\partial x_{j}}%
(s,N_{s})\Big(\frac{d}{ds}b^{j}_s ds+\sigma _{s}^{j}\delta X_{s}^{j}\Big) 
\notag \\
& -\frac{1}{2}\sum_{j=1}^{n}\int_{t}^{T}\frac{\partial ^{2}u}{\partial
x_{j}^{2}}(s,N_{s})\frac{d}{ds}Var(N_{s}^{j})ds.
\label{ItodiagonalBSDE}
\end{align}%
An application of $(\ref{PDE})$ to the second term on the right hand side of 
$(\ref{ItodiagonalBSDE})$ yields 
\begin{eqnarray}
u(t,N_{t}) &=&u(T,N_{T})-\int_{t}^{T}\Big(f(s)+A_{1}(s)u(s,N_{s})+%
\sum_{j=1}^{n}A_{2}^{j}(s)\sigma _{s}^{j}\frac{\partial u}{\partial x_{j}}%
(s,N_{s})\Big)ds  \notag \\
&&-\sum_{j=1}^{n}\int_{t}^{T}\sigma _{s}^{j}\frac{\partial u}{\partial x_{j}}%
(s,N_{s})\delta X_{s}^{j}.  \label{ItodiagonalBSDEII}
\end{eqnarray}%
We get $(\ref{BSDE})$ by setting $Y_{t}:=u(t,N_{t})$ and $Z_{t}^{j}:=-\sigma
_{t}^{j}\partial /\partial x_{j}u(t,N_{t})$, i.e. $(Y,(Z^{1},...,Z^{n}))$
solves $(\ref{BSDE})$ and is adapted to $\mathbb{F}.$ As in Section 1.2 we
consider this equation for $t\in \lbrack t_{0},T],$ for some fixed $%
t_{0}\geqslant 0$.\newline

In order to solve $(\ref{PDE})$ explicitely, we have to assume, in addition
to \emph{(H1) }and \emph{(H2)} for the kernels $K^{j},$ some regularity and
integrability conditions. \emph{(H5)} will be discussed later$.$ 

\emph{(H3)\qquad }There exist constants $c,C>0$ such that $c<\sigma ^{j}<C,$ 
$j=1,...,n,$ and $A_{1},f$, $A_{2}:=(A_{2}^{1},...,A_{2}^{n})$ are bounded.

\emph{(H4)\qquad }$g$ is continous, and there exist positive constants $%
c^{\prime }$ and $\lambda ^{\prime }<\underset{j=1,...,n}{\min }(16\underset{%
t\in \lbrack 0,T]}{\sup }Var(N_{t}^{j}))^{-1}$ such that $\mid
g(x)\mid \leqslant c^{\prime }e^{\lambda ^{\prime }\mid x\mid ^{2}}$ for all 
$x\in 
\mathbb{R}
^{n}.$

\emph{(H5)\qquad }$d /dtVar(N_{t}^{j})>0$ for all $%
t\in \lbrack t_{0},T]$ and $\int_{t_{0}}^{T}(Var(N_{T}^{j})-%
Var(N_{t}^{j}))^{-1/2}dt<\infty ,$ $j=1,...,n.$

\begin{thm}
Assume that (H1)-(H5) hold, and let $v(t,z)=(2\pi t)^{-1/2}\exp (-z^{2}/2t).$
Then $(Y,(Z^{1},...,Z^{n}))$ given by $(\ref{YetZ})$ solves $(\ref{BSDE}),$
where the classical solution of $(\ref{PDE})$ is given by%
\begin{align}
u(t,x) &=-\int_{t}^{T}\exp \Big(\int_{s}^{t}A_{1}(r)dr\Big)f(s)ds+\exp \Big(-%
\int_{t}^{T}A_{1}(s)ds\Big)  \notag \\
&\times \int_{\mathbb{R}^n} g(y)\prod\limits_{j=1}^{n}v\Big(Var(N_{T}^{j})-%
Var(N_{t}^{j}),x_{j}-\int_{t}^{T}(\sigma
_{s}^{j}A_{2}^{j}(s)-\frac{d}{ds}b^{j}_s)ds-y_{j}\Big)dy.  \label{solution}
\end{align}
\label{theosolutionEDP}
\end{thm}

An explicit calculation of $\partial /\partial tu(t,x)$ shows that $u$ is in
fact the (classical) solution of $(\ref{PDE}).$ We show now that $u$
verifies the growth condition $(\ref{F}).$

\bigskip

\begin{lem}
Let $u$ be given by $(\ref{solution})$. Then there are positive constants $M$
and $\lambda <\underset{j=1,...n}{\min }(4\underset{t\in \lbrack 0,T]}{\sup }%
Var(N_{t}^{j}))^{-1}$ such that
\begin{enumerate}
\item[1)] $\mid u(t, x)\mid \leqslant Me^{\lambda \mid x\mid ^{2}},$ $(t, x)\in \lbrack t_{0},T]\times \mathbb{R}^{n},$\newline

\item[2)] $\Big|\frac{\partial u}{\partial x_{j}}(t, x)\Big|\leqslant M(%
Var(N_{T}^{j})-Var(N_{t}^{j}))^{-1/2}e^{\lambda \mid
x\mid ^{2}},$ $(t, x)\in \lbrack t_{0},T)\times \mathbb{R}^{n},$\newline

\item[3)] $\Big|\frac{\partial ^{2}u}{\partial x_{j}^{2}}(t, x)\Big|\leqslant
M(Var(N_{T}^{j})-Var(N_{t}^{j}))^{-1}e^{\lambda \mid
x\mid ^{2}},$ $(t,x)\in \lbrack t_{0}, T)\times \mathbb{R}^{n}.$\\
\end{enumerate}
\end{lem}

Proof. We give the proof of 1), the proofs of 2) and 3) are similar. We
write $r_{s}^{j}$ for $\sigma _{s}^{j}A_{2}^{j}(s)-\frac{d}{ds}b^{j}_{s}$%
\ and $D_{t}^{j}$ for $Var(N_{T}^{j})-Var(N_{t}^{j})$
in $(\ref{solution}).$ Then, by the change of variables $z_{j}=x_{j}-%
\int_{t}^{T}r_{s}^{j}ds-y_{j},$ 
\begin{equation*}
\Big|u(t,x)\Big|\leqslant C\displaystyle\int_{\mathbb{R}^{n}}exp\Big(\lambda
^{\prime }\left\vert x-\int_{t}^{T}r_{s}ds-z\right\vert ^{2}\Big)%
\prod_{j=1}^{n}v(D_{t}^{j},z_{j})dy.
\end{equation*}%
Notice that 
\begin{equation}
\int_{\mathbb{R}^{n}}exp\Big(\lambda ^{\prime }\left\vert
x-\int_{t}^{T}r_{s}ds-z\right\vert ^{2}\Big)\leqslant exp\Big(2\lambda ^{\prime
}\left\vert x-\int_{t}^{T}r_{s}ds\right\vert ^{2}\Big)exp\Big(2\lambda
^{\prime }\mid z\mid ^{2}\Big)
\end{equation}

Therefore

\begin{align*}
\Big|u(t,x)\Big|& \leqslant Cexp\Big(2\lambda ^{\prime }\left\vert
x-\int_{t}^{T}r_{s}ds\right\vert ^{2}\Big)\int_{\mathbb{R}%
^{n}}\prod_{j=1}^{n}exp\Big(-\frac{\frac{1}{2}-2\lambda ^{\prime }Var%
(N_{T}^{j})}{D_{t}^{j}}z_{j}^{2}\Big)(2\pi D_{t}^{j})^{-1/2}dz \\
& =Cexp\Big(2\lambda ^{\prime }\left\vert x-\int_{t}^{T}r_{s}ds\right\vert
^{2}\Big)\prod_{j=1}^{n}\left( 1-4\lambda ^{\prime }Var%
(N_{T}^{j})\right) ^{-1/2} \\
& \leqslant Cexp\Big(2\lambda ^{\prime }\left\vert
x-\int_{t}^{T}r_{s}ds\right\vert ^{2}\Big)\leqslant Cexp\Big(4\lambda
^{\prime }\mid x\mid ^{2}\Big)exp\Big(4\lambda ^{\prime }\left\vert
\int_{t}^{T}r_{s}ds\right\vert ^{2}\Big),
\end{align*}%
and 1) follows by choosing $M=exp\Big(4\lambda ^{\prime }\vert%
\int_{t}^{T}r_{s}ds\vert^{2}\Big)$ and $\lambda =4\lambda ^{\prime }.$ \ \ $%
\blacksquare $\newline

\bigskip

Proof of Theorem \ref{theosolutionEDP} It remains to show that $(Y,Z)$ satisfies the BSDE $(%
\ref{BSDE}).$ By the preceding lemma $u$ satisfies (\ref{ItodiagonalBSDEII})
with $T$ replaced by $T-\varepsilon .$ Since $\mathbb{E}\exp (2\lambda \mid
N_{t}\mid ^{2})$ is bounded on $[0,T]$ for $\lambda <\underset{j=1,...n}{%
\min }(4\underset{t\in \lbrack 0,T]}{\sup }Var(N_{t}^{j}))^{-1},$
we obtain%
\begin{equation*}
\mathbb{E}\int_{t}^{T}\mid
f(s)+A_{1}(s)u(s,N_{s})+\sum_{j=1}^{n}A_{2}^{j}(s)\sigma _{s}^{j}\frac{%
\partial u}{\partial x_{j}}(s,N_{s})\mid ds\leqslant
\sum_{j=1}^{n}\int_{t}^{T}(Var(N_{T}^{j})-Var%
(N_{t}^{j}))^{-1/2}dt<\infty 
\end{equation*}%
for all $t<T$ by \emph{(H5). }By continuity of $u$ and $N,$ the terms in the
first line of (\ref{ItodiagonalBSDEII}) with $T$ replaced by $T-\varepsilon $
converge to the terms in (\ref{ItodiagonalBSDEII}) as $\varepsilon
\rightarrow 0.$ The divergence integrals converge too in the sense 
\begin{equation*}
\mathbb{E}\left[ \int_{t}^{T-\varepsilon }\sigma _{s}^{j}\frac{\partial u}{%
\partial x_{j}}(s,N_{s})\delta X_{s}^{j}F\right] \underset{\varepsilon
\rightarrow 0}{\rightarrow }\mathbb{E}\left[ \int_{t}^{T}\sigma _{s}^{j}%
\frac{\partial u}{\partial x_{j}}(s,N_{s})\delta X_{s}^{j}F\right] 
\end{equation*}%
for all $j=1,...,n$ and $F\in \mathcal{S}.$ \ \ $\blacksquare $

\textbf{Remark. }We discuss now the hypothesis \emph{(H5).} The positivity
of $\frac{d }{dt}Var(N_{t}^{j})$ means that $R(t):=%
Var(N_{t}^{j})$ is (strictly) increasing on $[t_{0},T]$. We note
that

\begin{equation*}
\frac{d }{dt}Var(N_{t}^{j})=\frac{d }{%
dt}\int_{0}^{t}(K_{t}^{\ast ,j}\sigma
^{j})_{s}^{2}ds=2\int_{0}^{t}\sigma _{t}^{j}\sigma _{u}^{j}\int_{0}^{u}\frac{%
\partial K^{j}}{\partial t}(t,s)\frac{\partial K^{j}}{\partial u}%
(u,s)dsdu=2\sigma _{t}^{j}\int_{0}^{t}\sigma _{u}^{j}\phi ^{j}(t,u)du.
\end{equation*}%
Since $\sigma ^{j}>0$ by \emph{(H3}$),$ a sufficient condition for $\frac{%
d}{dt}Var(N_{t}^{j})>0$ for all $t\in \lbrack
t_{0},T]$ is $\phi ^{j}>0$ on $[0,T]^{2}\setminus[0,t_{0}]^{2}.$ This is the case if 
$\frac{\partial K^{j} }{\partial u}(u,s)>0$ for all $(u,s)\in \lbrack
0,T]^{2},u>s.$ For mbm (Example 1) such cases can easily be found for Hurst
functions $h$ that are increasing or decreasing on $[0,T]$, but other
examples of mbm show that $R$ is increasing on $[t_{0},T]$ is possible for suitable
values of $t_{0}$ even if $h$ is only of bounded variation. We note that for
fractional brownian motion $B^{H}$ with Hurst index $H>1/2$%
\begin{equation*}
\phi ^{H}(t,u)=\frac{\partial ^{2}}{\partial t\partial u}\mathbb{E}%
B_{t}^{H}B_{u}^{H}=C_{H}\mid t-u\mid ^{2H-2}>0,
\end{equation*}

where $C_{H}>0$ is a constant depending on $H.$

Let us comment now on the hypothesis of integrability of ($%
D_{t}^{j})^{-1/2}:=(Var(N_{T}^{j})-Var%
(N_{t}^{j}))^{-1/2}$ near $T.$ We have

\begin{align*}
D_{t}^{j}& =\mathbb{E}\Big[\Big(\displaystyle\int_{0}^{T}\delta
W_{s}^{j}\int_{s}^{T}\sigma _{r}^{j}\frac{\partial K^{j}}{\partial r}(r,s)dr%
\Big)^{2}-\Big(\displaystyle\int_{0}^{t}\delta W_{s}^{j}\int_{s}^{t}\sigma
_{r}^{j}\frac{\partial K^{j}}{\partial r}(r,s)dr\Big)^{2}\Big] \\
& =\mathbb{E}\Big(\displaystyle\int_{0}^{T}\delta
W_{s}^{j}\int_{s}^{T}\sigma _{r}^{j}\frac{\partial K^{j}}{\partial r}(r,s)dr+%
\displaystyle\int_{0}^{t}\delta W_{s}^{j}\int_{s}^{t}\sigma _{r}^{j}\frac{%
\partial K^{j}}{\partial r}(r,s)dr\Big) \\
& \times \Big(\displaystyle\int_{0}^{T}\delta W_{s}^{j}\int_{s}^{T}\sigma
_{r}^{j}\frac{\partial K^{j}}{\partial r}(r,s)dr-\displaystyle%
\int_{0}^{t}\delta W_{s}^{j}\int_{s}^{t}\sigma _{r}^{j}\frac{\partial K^{j}}{%
\partial r}(r,s)dr\Big)\Big].
\end{align*}

An explicit calculation shows 
\begin{align*}
D_{t}^{j}& =\int_{0}^{t}dr\int_{t}^{T}dr^{\prime }\sigma _{r}^{j}\sigma
_{r^{\prime }}^{j}\phi ^{j}(r,r^{\prime })+\int_{0}^{t}ds\left(
\int_{t}^{T}dr\sigma _{r}^{j}\frac{\partial K^{j}}{\partial r}(r,s)\right)
^{2}+\int_{t}^{T}ds\left( \int_{s}^{T}dr\sigma _{r}^{j}\frac{\partial K^{j}}{%
\partial r}(r,s)\right) ^{2} \\
& =:A_{t}^{1}+A_{t}^{2}+A_{t}^{3}.
\end{align*}

Under the hypothesis \emph{(H3) }for $\sigma ^{j}$ and if $\phi ^{j}>0,$ a
sufficient condition for $\int_{t_{0}}^{T}$ ($D_{t}^{j})^{-1/2}dt<\infty $ is%
\begin{equation*}
A_{t}^{3}=\int_{t}^{T}(K_{T}^{\ast ,j}\sigma ^{j})_{s}^{2}ds\geqslant
c(T-t)^{a}
\end{equation*}

for some constant $c>0$ and $a\in (0,2)$ as $t\nearrow T.$ For fractional
brownian motion $B^{H}$ with $H>1/2$ this condition is satisfied with $%
a=H+1/2$. For the Volterra processes in Examples 1-3 this condition is
satisfied with $a=2h(T)$ if $h$ is such that $\frac{\partial K}{\partial u}%
(u,s)>0,$ for $(u,s)\in (t_{0},T)^{2},u>s$.

\section*{Acknowledgments}

The authors thank the program Hubert Curien "Utique" of the 'French Ministry
of Foreign Affairs' and the 'Tunisian Ministry of Education and Research'
for the financial support. \newline

\bibliographystyle{plain}
\bibliography{KnaniDozzibiblio}

\end{document}